\documentclass[12pt]{amsart}
\usepackage{amssymb}
\usepackage{array}
\usepackage{latexsym}
\usepackage{enumerate}
\usepackage{amsmath}
\usepackage{amsfonts}
\usepackage{amsthm}
\usepackage{verbatim}
\usepackage[english]{babel}
\newtheorem{theorem}{Theorem}[section]
\newtheorem{proposition}[theorem]{Proposition}
\newtheorem{lemma}[theorem]{Lemma}

{\theoremstyle{definition}
\newtheorem*{definition*}{Definition}

\newtheorem{rem}[theorem]{Remark}
\newtheorem*{proposition*}{Proposition}
\newtheorem*{corollary*}{Corollary}
\newtheorem*{lemma*}{Lemma}

\def\cC{\mathcal C}

\def\cE{\mathcal E}
\def\cF{\mathcal F}

\def\cL{\mathcal L}

\def\cU{\mathcal U}

\def\cX{\mathcal X}
\def\cY{\mathcal Y}
\def\cZ{\mathcal Z}

\def\Aut{\mbox{\rm Aut}}
\def\K{\mathbb{K}}
\def\PG{{\rm{PG}}}
\def\ord{\mbox{\rm ord}}
\def\deg{\mbox{\rm deg}}

\def\Aut{\mbox{\rm Aut}}

\newcommand{\PSL}{\mbox{\rm PSL}}
\newcommand{\PGL}{\mbox{\rm PGL}}
\newcommand{\AGL}{\mbox{\rm AGL}}
\newcommand{\PSU}{\mbox{\rm PSU}}
\newcommand{\PGU}{\mbox{\rm PGU}}

\newcommand{\Sz}{\mbox{\rm Sz}}

\newcommand{\SU}{\mbox{\rm SU}}

\newcommand{\aut}{\mbox{\rm Aut}}

\newcommand{\ga}{\alpha}
\newcommand{\gb}{\beta}
\newcommand{\g}{\gamma}

\newcommand{\gep}{\epsilon}

\newcommand{\gl}{\lambda}

\newcommand{\gz}{\zeta}

\newcommand{\gO}{\Omega}

\newcommand{\ha}{{\textstyle\frac{1}{2}}}
\newcommand{\thi}{\textstyle\frac{1}{3}}

\newcommand{\tth}{\textstyle\frac{2}{3}}

\newcommand{\bv}{{\bf v}}

\sloppy

\begin{document}
\def\amshead{
\title[Automorphisms of zero $p$-rank curves]{}
\author{M.~Giulietti and G.~Korchm\'aros}}

 \thanks{Research supported by  the Italian
    Ministry MURST, Strutture geometriche, combinatoria e loro
    applicazioni, PRIN 2006-2007}

\thanks{2000 {\em Math. Subj. Class.}: 14H37}

\begin{center}
\begin{LARGE}
\bf{
Automorphism groups of algebraic\\ \vspace*{.18cm} curves
with $p$-rank zero}\\
\end{LARGE}
\end{center}
\amshead

\vspace*{-.5cm}
    \begin{abstract} In positive characteristic, algebraic curves can have many more automorphisms
    than expected from the classical Hurwitz's bound. There even exist algebraic curves of
    arbitrary high genus $g$ with more than $16g^4$ automorphisms. It has been observed on many occasions that the most anomalous
    examples invariably have zero $p$-rank. In this paper, the $\K$-automorphism group $\aut(\cX)$ of a zero $2$-rank algebraic curve
    $\cX$ defined over an algebraically closed field $\K$ of characteristic $2$ is investigated. The main result is that if the curve has genus $g\geq 2$
    and $|\aut(\cX)|>24g^2$, then $\aut(\cX)$ has a fixed point on $\cX$, apart from few
    exceptions. In the exceptional cases the possibilities for $\aut(\cX)$ and $g$ are determined.
    \end{abstract}
\maketitle

    \section{Introduction}\label{sec1}
Let $\cX$ be a (projective, geometrically irreducible,
non-singular) algebraic curve defined over an algebraically closed
groundfield $\K$ of characteristic $p\geq 0.$ Let $\K(\cX)$ be the
field of rational functions (the function field of one variable
over $\K$) of $\cX$. The $\K$--automorphism group $\aut(\cX)$ of
$\cX$ is defined to be the automorphism group $\aut(\K(\cX))$
consisting of those automorphisms of $\K(\cX)$ which fix each
element of $\K$.

By a classical result, $\aut(\cX)$ is finite if the genus $g$ of
$\cX$ is at least two, see
\cite{schmid1938,iwasawatamagawa1951,iwasawatamagawa1951b,iwasawatamagawa1951c,roquette1952}.

It is known that every finite group occurs in this way, since for
any groundfield $\K$ and any finite group $G$, there exists $\cX$
such that $\aut(\cX)\cong G$, see Greenberg \cite{greenberg1974}
for $\K={\mathbb{C}}$ and Madden and Valentini
\cite{maddenevalentini1983} for $p\geq 0$, see also Madan and
Rosen  \cite{madanrosen1992} and Stichtenoth
\cite{stichtenoth1984}.




This raises a general problem for groups and curves:
Determine the finite groups that can be realized as the
$\K$-automorphism group of some curve with a given invariant, such
as genus $g$, $p$-rank (for $p>0$) or number of
${\mathbb{F}}_q$-rational places (for $\K=\bar{\mathbb{F}}_q$),
where ${\mathbb{F}}_q$ stands for the finite field of order $q$.


In the present paper we deal with zero $p$-rank curves and their
$\K$--automorphism groups. Interestingly, zero $p$-rank curves
with large $\K$--automorphism groups are known to exist, see
Roquette \cite{roquette1970}, Stichtenoth
\cite{stichtenoth1973I,stichtenoth1973II}, Henn \cite{henn1978},
Hansen and Petersen \cite{hansen1993}, Garcia, Stichtenoth and
Xing \cite{garcia-stichtenoth-xing2001}, \c{C}ak\c{c}ak and
\"Ozbudak \cite{cakcak-ozbudak2004}, Giulietti, Korchm\'aros and
Torres \cite{giulietti-korchmaros-torres2004}, Lehr and Matignon
\cite{lehr-matignon2005}, Giulietti and Korchm\'aros
\cite{segovia}. Those for $p=2$ are the non-singular models of the
plane curves listed below.
\begin{enumerate}
    \item[\rm(I)] $\bv(Y^2+Y+X^{2^k+1})$, a hyperelliptic curve of genus
$g=2^{k-1}$ with $\aut(\cX)$ fixing a point of $\cX$ and
$|\Aut(\cX)|=2^{2k+1} (2^k+1);$
    \item[\rm(II)] $\bv(Y^n+Y - X^{n+1})$, the Hermitian curve of degree $n+1$ with $\aut(\cX)\cong
    \PGU(3,n),$
    $n$ a power of $2$;
    \item[\rm(III)] $\bv(X^{n_0}(X^n+X)+Y^n+Y)$, the DLS curve (Deligne-Lusztig curve of Suzuki type) of
genus $g=n_0(n-1)$ with $\aut(\cX)\cong \Sz(n)$ where $\Sz(n)$ is
the Suzuki group, $n_0=2^r, r\geq 1, n=2n_0^2$.

\item[\rm(IV)] $\bv(Y^{n^3+1}+(X^n+X)(\sum_{i=0}^n
(-1)^{i+1}X^{i(n-1)})^{n+1})$, a curve of genus
$g=\ha\,(n^3+1)(n^2-2)+1$ with $\aut(\cX)$ containing a subgroup
isomorphic to $\SU(3,n)$, $n$ a power of $2$.
\end{enumerate}

The essential idea in our investigation is to deduce from the zero
$p$-rank condition the purely group theoretic property that every
$\K$--automorphism of $\cX$ of order $p$ has exactly one fixed
point. Then deeper results from Group Theory can be used to
determine the structure and action of $\aut(\cX)$ on the set of
points of $\cX$. This idea works well for $p=2$
and produces the following result.
\begin{theorem}
\label{th1main} Let $\cX$ be a zero $2$-rank algebraic curve of
genus $g\geq 2$ defined over an algebraically closed groundfield
$\K$ of characteristic $2$.  Let $G$ be a subgroup of $\aut(\cX)$
of even order.
Let $S$ be the subgroup generated by all elements of $G$ of order
a power of $2$. If no point of $\cX$ is fixed by $G$, then one of
the following holds$:$
\begin{enumerate}
\item[\rm(a)] $S$ is isomorphic to one of the groups below$:$
\begin{equation}
\label{lingp}
 \PSL(2,n),\, \PSU(3,n),\, {\mbox{\rm SU}}(3,n),\,\Sz(n)\,\,  \mbox{with $n =2^r\ge 4$};
\end{equation}
\item[\rm(b)] $S=O(S)\rtimes S_2$, where $O(S)$ is the largest normal subgroup of odd
order of $S$, and $S_2$ is a Sylow $2$-subgroup which is either a cyclic group or a generalised quaternion group.
\end{enumerate}
If $G$ fixes a point of $\cX$ then $G=S_2\rtimes H,$ with a Sylow
subgroup $S_2$ of $G$.
\end{theorem}
Theorem \ref{th1main} is the essential tool in the present
investigation on the spectrum of the genera of zero $2$-rank
curves with large automorphism groups. Our results are summarized
below.

Let $\cX$ be a zero $2$-rank curve. If $\aut(\cX)$ is a solvable
group and has no fixed point on $\cX$, then $|\aut(\cX)|\leq
24g^2$, see Theorem \ref{th2hering}.

Assume that $\aut(\cX)$ is non-solvable. Then $S$ coincides with
the commutator group $\aut(\cX)'$ of $\aut(\cX)$, see Theorem
\ref{th2main}. The possible genera of $\cX$ can be computed from
the order of $\aut(\cX)'$ provided that $\aut(\cX)$ is big enough,
namely whenever
both
\begin{equation}
\label{eqfinal} {\mbox{\em{$|\aut(\cX)_P|$ is even and bigger than
$3g$ for at least one point $P\in \cX$}}}
\end{equation}
and
\begin{equation}
\label{eqfinalbis} |\aut(\cX)'|\geq 6(g-1)
\end{equation}
hold,
see Theorem \ref{thmain}. Such possibilities for $g$ are the
following:
\begin{itemize}
    \item[\rm(i)]\quad $\aut(\cX)'=\PSL(2,n)$ and $g=\ha(t-1)(n-1)$ with $t|(n+1).$
    \item[\rm(ii)]\quad $\aut(\cX)'=\PSU(3,n)$ and either
            $$g=\ha\,(n-1)(t(n+1)^2-(n^2+n+1))$$ with
    $t|(n^2-n+1)/\mu$, and $\mu={\rm{gcd}}(3,n+1)$, or
$$
g=\ha (n-1)\left(\frac{t(n^3+1)}{\mu}-(n^2+n+1)\right),
$$
with $t\mid n+1$;
    in the former case,   $t=1$ only occurs when $\cX$ is as in (II).
    \item[\rm(iii)]
\quad $\aut(\cX)'=\Sz(n)$ and
    either $$g\,=\ha\,[(t-1)(n^2-1)-2tn_0(n-1)]$$ with $t|(n
    +2n_0+1)$, or $$g\,=\ha\,[(t-1)(n^2-1)+2tn_0(n-1)]$$ with $t|(n
    -2n_0+1)$; in
     the latter case $t=1$ only occurs when $\cX$ is as in (III).

    \item[\rm(iv)]
   \quad $\aut(\cX)'=\SU(3,n)$ with $3|(n+1)$ and either
$$g=\ha (n-1)[3t(n+1)^2-(n^2+n+1)]$$ with $t|(n^2-n+1)/3$, or
$$
g=\ha (n-1)\left(t(n^3+1)-(n^2+n+1)\right),
$$
with $t\mid (n+1)$; in the former case $t=n^2-n+1$ (equivalently,
in latter case $t=n+1$) occurs when $\cX$ is as in (IV).
\end{itemize}
Both technical conditions (\ref{eqfinal}) and (\ref{eqfinalbis})
are
satisfied when $\aut(\cX)$ is large enough, see
Lemma \ref{global}. Therefore, the above results have the
following corollary.
\begin{theorem}
\label{th1maincor} Let $\cX$ be a zero $2$-rank algebraic curve of
genus $g\geq 2$  defined over an algebraically closed groundfield
$\K$ of characteristic $2$. Assume that $\aut(\cX)$ fixes no point
on $\cX$. If $|\aut(\cX)|>24g^2$, then one of the above four cases
{\rm (i)-(iv)} occurs, and  $\aut(\cX)$ contains a cyclic normal
subgroup
$N$
of odd order such that the factor group $\aut(\cX)/N$ is
isomorphic to one of the groups $\PSL(2,n)$, $\PSU(3,n)$,
$\PGU(3,n)$ and $\Sz(n)$.
\end{theorem}
For $\PSL(2,n)$ and $\Sz(n)$, the above central extension is
splitting, see Theorem \ref{thmaincor1}. In any case, $N$
coincides with the subgroup of $\aut(\cX)$ fixing point-wise the
unique non-tame orbit of $\aut(\cX)$ on $\cX$.

It may be asked whether similar results for zero $p$-rank curves
can hold for odd $p$. Here we limit ourselves to extend the bound
$|G|\leq 24g^2$ to solvable $\K$--automorphism groups $|G|$ whose
order is divisible for $p^2$.

Surveys on $\K$-automorphism groups of curves in positive
characteristic are found in \cite{nakajima2002} and
\cite{hirschfeld-korchmaros-torres2008}.

    \section{Background}\label{sec2}

$\aut(\cX)$ has a faithful permutation representation on the set
of all points of $\cX$ (equivalently on the set of all places of
$\K(\cX)$). The orbit
$$o(P)=\{Q\mid Q=P^\alpha, \alpha\in G\}$$
is {\em long} if $|o(P)|=|G|$, otherwise $o(P)$ is {\em short} and
$G_P$ is non-trivial.

If $G$ is a finite subgroup of $\aut(\cX)$, the subfield
$\K(\cX)^G$ consisting of all elements of $\K(\cX)$ fixed by every
element in $G$, is also a function field of one variable over
$\K$. Let $\cY$ be a non-singular model of $\K(\cX)^G$, that is, a
(projective, geometrically irreducible, non-singular) algebraic
curve defined over $\K$ with function field $\K(\cX)^G$. Usually,
$\cY$ is called the quotient curve of $\cX$ by $G$ and denoted by
$\cX/G$. The covering $\cX\mapsto \cY$ has degree $|G|$ and the
field extension $\K(\cX)/\K(\cX)^G$ is Galois.

If $P$ is a point of $\cX$, the stabiliser $G_P$ of $P$ in $G$ is
the subgroup of $G$ consisting of all elements fixing $P$. For
$i=0,1,\ldots$, the $i$-th ramification group $G_P^{(i)}$ of $\cX$
at $P$ is
$$G_P^{(i)}=\{\alpha\mid \ord_P(\alpha(t)-t)\geq i+1, \alpha\in
G_P\}, $$ where $t$ is a uniformizing element of $\cX$ at $P$.
Here $G_P^{(0)}=G_P$ and $G_P^{(1)}$ is the unique Sylow
$p$-subgroup of $G_P$. If $\cX/G_P^{(1)}$ is non-rational, then
$|G_P^{(1)}|\leq g$, by a result due to Stichtenoth
\cite{stichtenoth1973I}. Also, $G_P=G_P^{(1)}\rtimes H$ with $H$ a
cyclic group whose order is prime to $p$. Furthermore, for $i\geq
1$, $G_P^{(i)}$ is a normal subgroup of $G_P$ and the factor group
$G_P^{(i)}/G_P^{(i+1)}$ is an elementary abelian $p$-group. For
$i$ big enough, $G_P^{(i)}$ is trivial.

For any point $Q$ of $\cX$, let $e_Q=|G_Q|$ and
$$d_Q=\sum_{i\geq 0}(|G_Q^{(i)}|-1).$$ Then $d_Q\geq e_Q-1$
and equality holds if and only if $\gcd(p,|G_Q|)=1.$

Let $g'$ be the genus of $\cX/G$. From the Hurwitz genus formula,
    \begin{equation}
    \label{eq1}
2g-2=|G|(2g'-2)+\sum_{Q\in \cX} d_Q.
    \end{equation}

If $G$ is tame, that is $p\nmid |G|$, or more generally for $G$
with $p\nmid e_Q$ for every $Q\in\cX$, Equation (\ref{eq1}) is
simpler and may be written as
    \begin{equation}
    \label{eq2}
2g-2=|G|(2g'-2)+\sum_{i=1}^k (|G|-\ell_i)
    \end{equation}
where $\ell_1,\ldots,\ell_k$ are the sizes of the short orbits of $G$ on $\cX$.

Let $G_P=G_P^{(1)}\rtimes H$. The following upper bound on $|H|$
depending on $g$ is due to Stichtenoth \cite{stichtenoth1973I}:
\begin{equation*}
\label{stichtenotharchivsatz2} |H|\leq 4g+2.
\end{equation*}

For any abelian subgroup $G$ of $\aut(\cX)$, Nakajima
\cite{nakajima1987bis} proved that
\begin{equation*}
|G|\leq \left \{
\begin{array}{ll}
4g+4 & \mbox{for\quad  $p\neq 2,$} \\
4g+2 & \mbox{for\quad  $p=2.$}
\end{array}
\right.
\end{equation*}

An equation similar to (\ref{eq2}) holds for any $p$-group $S$ of
$\K$--automorphisms of $\cX$ whenever the genus $g$ is replaced by
the $p$-rank. Recall that the $p$-rank $\g$ of a curve is its
Hasse-Witt invariant and that $0\leq \g\leq g$. The elementary
abelian subgroup ${\rm{Pic}}_0(\cX,p)$ of elements of order $p$ in
the zero divisor class group ${\rm{Pic}}_0(\cX)$ has dimension
$\g$, that is, ${\rm{Pic}}_0(\cX,p)\cong (\mathbb{Z}_p)^\g$. If
$\g'$ denotes the $p$-rank of $\cX/S$, the Deuring-Shafarevich
formula is
\begin{equation}
    \label{eqds}
\g-1=|
S
|(\g'-1)+\sum_{i=1}^k (|S|-\ell_i)
    \end{equation}
where $\ell_1,\ldots,\ell_k$ are the short orbits of $S$ on $\cX$.
Upper bounds on $|S|$ depending on $\g$ are found in
\cite{nakajima1987,nakajima1987bis}.

Let $m$ be the smallest non-gap at $P$, and suppose that $\xi\in
\K(\cX)$ has a pole of order $m$ at $P$ and is regular elsewhere. If
$g\in G_P$, then $g(\xi)=c\xi+d$ with $c,d\in \K$ and $c\neq 0$.

Let $\cL$ be the projective line over $\K$. Then $\aut(\cL)\cong
\PGL(2,\K),$ and $\aut(\cL)$ acts on the set of all points of
$\cL$ as $\PGL(2,\K)$ naturally on $\PG(1,\K)$. In particular, the
identity of $\aut(\cL)$ is the only $\K$--automorphism in
$\aut(\cL)$ fixing at least three points of $\cL$. Every
$\K$--automorphism $\ga\in\aut(\cL)$ fixes a point; more
precisely, $\ga$ fixes either one or two points according as its
order is $p$ or relatively prime to $p$. Also, $G_P^{(1)}$ is an
infinite elementary abelian $p$-group. For a classification of
subgroups of $\PGL(2,\K)$, see \cite{maddenevalentini1982}.

Let $\cE$ be an elliptic curve. Then $\aut(\cE)$ is infinite;
however for any point $P\in \cE$ the stabiliser of $P$ is rather
small, namely
\begin{equation*}
|\aut(\cE)_{P}| = \left \{
\begin{array}{ll}
 2,4,6 & \mbox{\quad when $p\neq 2,3,$} \\
 2,4,6,12 & \mbox{\quad when $p=3,$} \\
 2,4,6,8,12,24 & \mbox{\quad when $p=2.$}
\end{array}
\right.
\end{equation*}

Let $\cF$ be a (hyperelliptic) curve of genus $2$
having more than one Weierstrass point.
For any solvable subgroup $G$ of $\aut(\cF)$, Nakajima's bound
together with some elementary facts on finite permutation groups,
yield $|G|\leq 48$.

{}From finite group theory, the following results play a role in
the proofs.

{\em Huppert's classification theorem}, see \cite{huppertblackburn3-1982}, Chapter
XII, Section 7: Let $G$ be a solvable
$2$-transitive permutation group of degree $n$. Then $n$ is a power of some prime $p$,  and $G$ is a subgroup of the affine semi-linear
group ${\mathrm{A}}\Gamma {\mathrm{L}}(1,n)$, except possibly when
$n$ is $3^2,5^2,7^2,11^2,23^2$ or $3^4$.

{\em The Kantor-O'Nan-Seitz theorem}, see
\cite{kantor-o'nan-seitz1972}: Let $G$ be a $2$-transitive
permutation group of odd degree. If the $2$-point stabiliser of
$G$ is cyclic, then $G$ has either an elementary abelian regular
normal subgroup$,$ or $G$ is one of the following groups in
its natural $2$-transitive permutation representation$:$
$$\PSL(2,n),\,\PSU(3,n),\,\PGU(3,n),\Sz(n),$$
where $n\geq 4$ is a power of $2$.

{\em The natural $2$-transitive permutation representations of the
above linear groups:}
\begin{itemize}
\item[(A)] $G=\PSL(2,n)$, is the $\K$--automorphism group of
$\PG(1,n)$; equivalently, $G$ acts on the set $\Omega$ of all
${\mathbb{F}}_{n}$-rational points of the projective line defined
over ${\mathbb{F}}_n$.

\item[(B)] $G=\PGU(3,n)$ is the linear collineation group
preserving the classical unital in the projective plane
$\PG(2,n^2)$, see \cite{hirschfeld1998}; equivalently $G$ is the
$\K$--automorphism group of the Hermitian curve regarded as a
plane non-singular curve defined over the finite field
${\mathbb{F}}_n$ acting on the set $\Omega$ of all
${\mathbb{F}}_{n^2}$-rational points,
see (II). For $n\equiv 0,1 \pmod 3$, $\PSU(3,n)=\PGU(3,n)$. For
$n\equiv -1 \pmod 3$, $\PSU(3,n)$ is a subgroup of $\PGU(3,n)$ of
index $3$, and this is the natural $2$-transitive representation
of $\PSU(3,n)$.

\item[(C)] $G=\Sz(n)$ with $n=2n_0^2$, $n_0=2^r$ and $r\geq 1$, is
the linear collineation group of $\PG(3,n)$ preserving the Tits
ovoid, see \cite{tits1960,tits1962,hirschfeld1985}; equivalently
$G$ is the $\K$--automorphism group of the {\rm{DLS}} curve
regarded as a non-singular curve defined over the finite field
${\mathbb{F}}_n$ acting on the set $\Omega$ of all
${\mathbb{F}}_{n}$-rational points, see (III).
\end{itemize}

For each of the above linear groups, the structure of the
$1$-point stabilizer and its action in the natural $2$-transitive
permutation representation, as well as its $\K$--automorphism
group, are explicitly given in the papers quoted. In this paper,
we will need the following corollary.
\begin{lemma}
\label{stabi}
Let $n$ be a power of $2$ and $n\geq 4$.
Let $L$ be any of the  groups $\PSL(2,n)$, $\PSU(3,n)$, $\Sz(n)$
given in its natural $2$-transitive permutation representation on
a set $\Omega$. Let ${\bar G}$ be a subgroup of $\aut(L)$
containing $L$ properly, and choose a point $P\in\Omega$. If the
Sylow
$2$-subgroup
of ${\bar G}_P$ properly contains the Sylow $2$-subgroup of $L_P$,
then some non-trivial element of order $2$ contained in ${\bar
G}_P\setminus L_P$ fixes a point other than $P$. Otherwise, either
${\bar G}_P$ contains a non-cyclic subgroup of odd order, or
$L=\PSU(3,n)$ and ${\bar G}=\PGU(3,n)$ with $3|(n+1)$.
\end{lemma}

{\em Cyclic fix-point-free subgroups of some $2$-transitive
groups.} The following lemma is a corollary of the classification
of subgroups of $\PSU(3,n)$ and $\Sz(n)$.
\begin{lemma}
\label{backgr} Let $G$ be a $2$-transitive permutation group of
degree $n$. Let $U$ be a cyclic subgroup of $G$ which contains no
non-trivial element fixing a point.
    \begin{itemize}
    \item[{\rm (i)}] If $G=\PSU(3,n)$ in its natural $2$-transitive
    permutation representation, then $|U|$ divides either $n+1$ or
    $(n^2-n+1)/\mu$, with $\mu={\rm{gcd}}(3,n+1)$.
    \item[{\rm (ii)}] If $G=\Sz(n)$ in its natural $2$-transitive
    permutation representation, then $|U|$ divides either $n-2n_0+1$ or $n+2n_0+1$.
    \end{itemize}
\end{lemma}

\section{$\K$--automorphism groups of curves with {$\lowercase{p}$}--rank zero}
\label{sectionstructure}

In this section, subgroups of $\aut (\cX)$ whose order is
divisible by $p$ and in which
\begin{equation}
\label{condstar} \mbox{\em every element of order $p$ has exactly
one fixed point}
\end{equation}
are investigated. In terms of ramification, (\ref{condstar}) means
that for every $\ga\in G$ of order $p$ the covering $\cX\mapsto
\cX/\langle\ga \rangle$ ramifies at exactly one point.

\begin{rem}
\label{condstar2}  A suitable power $\gb$ of any non-trivial
element $\ga$ of a $p$-group has order $p$. Also, every fixed
point of $\ga$ is fixed by $\gb$ as well. Since $G_P^{(1)}$ is a
$p$-group, it is a subgroup of a Sylow $p$-subgroup $S_p$ of $G$.
Choose a non-trivial element $\gz$ from the centre of $S_p$.
 Then $\gz$ commutes with a non-trivial element
of $G_P^{(1)}$. This together with (\ref{condstar}) imply that
$\gz$ fixes $P$. Therefore, $\gz\in G_P^{(1)}$. In particular,
$\gz$ fixes no point of $\cX$ distinct from $P$. Let $\ga\in S_p$.
Then $\gz \ga=\ga \gz$ implies that
$$(P^\ga)^\gz=P^{\ga\gz}=P^{\gz\ga}=(P^\gz)^\ga=P^\ga$$ whence $P^\ga=P$. This shows that
every element of $S_p$ must fix $P$, and hence $S_p=G_P^{(1)}$.
Since the Sylow $p$-subgroups of $G$ are conjugate under $G$,
every $p$-element fixes exactly one point of $\cX$.

Therefore, condition (\ref{condstar}) is satisfied by a subgroup
$G$ of $\aut (\cX)$ if and only if every $p$-element in $G$ has
exactly one fixed point. Further, the non-trivial elements of a
Sylow $p$-subgroup $S_p$ of a $\K$-automorphism $G$ of $\aut
(\cX)$ satisfying (\ref{condstar}) have the same fixed point.
\end{rem}

Property (\ref{condstar}) is known to occur in several
circumstances, for instance in \cite{stichtenoth1973I}, Satz 1,
see also \cite{lehr-matignon2005} and
\cite{giulietti-korchmaros2007}. Another such circumstance is
described in the following two lemmas.

\begin{lemma}
\label{deuringshafarevic2} If $\cX$ has $p$-rank $0,$ then
{\rm (\ref{condstar})} holds in $\aut (\cX)$.
\end{lemma}
\begin{proof} Let $\ga\in \aut(\cX)$ have order $p$. Applying the
Deuring--Shafarevich formula (\ref{eqds}) to the group generated
by $\ga$ gives $-1=p(\g'-1)+m(p-1)$, where $\g'$ is the $p$-rank
of $\K(\cX)^{\ga}$ and $m$ is the number of fixed points of $\ga$.
This is only possible when $\g'=0$ and $m=1$.
\end{proof}

The converse of Lemma \ref{deuringshafarevic2} is not true in
general, a counterexample for $p=2$ being the non-singular model
of the irreducible curve
$$\cF:\quad \bv(Y^6+Y^5+Y^4+Y^3+Y^2+Y+1+X^3(Y^2+Y)).$$  The $2$-rank of such a curve $\cX$ is
equal to $4$, while the birational transformation $x'=x,\,y'=y+1$
is a $\K$--automorphism of $\K(\cX)$ which has only one fixed
point, namely that arising from the unique branch of $\cF$ tangent
to the infinite line.
A partial converse of Lemma \ref{deuringshafarevic2} is the
following result.

\begin{lemma}
\label{deuringshafarevic2bis} Assume that $\aut (\cX)$ contains a
$p$-subgroup $G$. If $\cY=\cX/G$ has $p$-rank
zero$,$ and {\rm(\ref{condstar})}  holds$,$ then $\cX$ has $p$-rank
zero, as well.
\end{lemma}

\begin{proof} Let $P$ be the fixed point of $G$.
Since $\g'=0$,  (\ref{eqds}) applied to $G$ gives
$\g-1=-|G|+|G|-1$, whence $\g=0$.
\end{proof}

Suppose that a subgroup $G$ of $\aut (\cX)$ has a Sylow
$p$-subgroup with property (\ref{condstar}). Then (\ref{condstar})
is satisfied by all Sylow $p$-subgroups of $G$. Suppose further
that $G$ has no normal Sylow $p$-subgroup. Let $S_p$ and $S_p'$
two distinct Sylow $p$-subgroups of $G$. By Remark
\ref{condstar2}, $S_p$ and $S'_p$ each have exactly one fixed
point, say $P$ and $P'$. Also, $S_p$ is the unique Sylow
$p$-subgroup of the stabiliser of $P$ under $G$, and
$G_{P}=S_p\rtimes H$. Further, $P\neq P'$; for, if $P=P'$,
then $S_p$ and $S_p'$ are two distinct Sylow $p$-subgroups in
$G_{P}$. However, this is impossible as $S_p$ is a normal
subgroup of $G_P$. Therefore, any two distinct Sylow $p$-subgroups
in $G$ have trivial intersection.

A subgroup $V$ of a group $G$ is a {\em trivial intersection set}
if, for every $g\in G$, either $V=g^{-1}Vg$ or $V\cap
g^{-1}Vg=\{1\}$.
A Sylow $p$-subgroup $T$ of $G$ is a trivial intersection set in
$G$ if and only if that $T$ meets every other Sylow $p$-subgroup
of $G$ trivially. Hence, Lemma \ref{deuringshafarevic2} has the
following corollary.

\begin{theorem}
\label{deuringshafarevic2tris} If $\cX$ has $p$-rank $0,$ then
every Sylow $p$-subgroup $S_p$ in $\aut (\cX)$ is a trivial
intersection set.
\end{theorem}

A classical result on trivial intersection sets is the following.
\begin{theorem}
\label{gow} {\rm (Burnside)} If a Sylow $p$-subgroup $S_p$ of a
finite solvable group is a trivial intersection set$,$ then either
$S_p$ is normal or cyclic$,$ or $p=2$ and $S_2$ is a generalised
quaternion group.
\end{theorem}
For a proof, see \cite{gow1975}.

Finite groups whose Sylow $2$-subgroups are trivial intersection
sets have been classified. This has been refined to groups
containing a subgroup of even order which intersects each of its
distinct conjugates trivially.

\begin{theorem}
\label{hering} {\rm (Hering \cite{hering1972})}  Let $V$ be a
subgroup of a finite group $G$ with trivial normaliser
intersection$;$ that is$,$ the normaliser $N_G(V)$ of $V$ in $G$
has the following two properties$:$
\begin{enumerate}
\item[\rm(a)] $V\cap x^{-1}Vx=\{1\}$ for all $x\in G\setminus
N_G(V);$ \item[\rm(b)] $N_G(V)\neq G$.
\end{enumerate}
If $|V|$ is even and $S$ is the normal closure of $V$ in $G$, then
the following hold.
    \begin{enumerate}
    \item[\rm(i)] $S=O(S)\rtimes V,$ the semi-direct
product of $O(S)$ by $V,$ with $V$ a Frobenius complement$,$
unless $S$ is isomorphic to one of the groups$,$
\begin{equation}
\label{eqhering}
 \PSL(2,n),\, \Sz(n),\, \PSU(3,n),\,{\mbox{\rm SU}}(3,n),
 \end{equation}
where $n$ is a power of $2$.
    \item[\rm(ii)] Let $\bar{G}$ be the
permutation group induced by $G$ on the set $\Omega$ of all
conjugates of $V$ under $G$. Then
    \begin{enumerate}
    \item[\rm ii(a)] $C_G(S)$ is the kernel of this representation$;$
    \item[\rm ii(b)]  $S$ is transitive on $\Omega;$ \item[\rm(c)]
$|\Omega|$ is odd$;$
    \item[\rm ii(d)] in the exceptional cases$,$
    \begin{enumerate}
    \item[\rm (1)] $|\Omega|=n+1,n^2+1,n^3+1,n^3+1;$
    \item[\rm (2)] $\bar{G}=G/C_G(S)$ contains a normal subgroup $L$ isomorphic to one of
the  groups $\PSL(2,n)$, $\Sz(n)$, $\PSU(3,n)$,
acting in its natural $2$-transitive permutation representation;
    \item[\rm (3)] $\bar{G}$
is isomorphic to an automorphism group of $L$ containing $L$.
    \end{enumerate}
    \end{enumerate}
    \end{enumerate}
\end{theorem}

Hering's result does not extend to subgroups of odd order with
trivial normal intersection. A major result which can be viewed as
a generalisation of Theorem \ref{gow} is found in
\cite{zhang1991}.

\begin{theorem}
\label{zhang1991} Let  $S_d,$ with $d>11,$ be a Sylow $d$-subgroup
of a finite group $G$ that is a trivial intersection set but not a
normal subgroup. Then $S_d$ is cyclic if and only if $G$ has no
composition factors isomorphic to either $\PSL(2,d^n)$ with $n>1$
or $\PSU(3,d^m)$ with $m\geq 1$.
\end{theorem}

Now, some consequences of the above results are stated.

Let
$p=2$. If (\ref{condstar})  is satisfied by a subgroup $G$ of
$\aut(\cX)$, then Theorem \ref{hering} applies to $G$ where $V$ is
a Sylow $2$-subgroup $S_2$ of $G$ (or any non-trivial subgroup of
$S_2)$. In fact, (a) is true, while (b) holds provided that $G$ is
larger than $G_{P}$, where $P$ is the fixed point of $S_2$.

When $G_P=G$, then $G=G_P^{(1)}\rtimes H$ with $p\nmid |H|$.
Otherwise, $G$ has no fixed point, and Hering's theorem determines
the abstract structure of the normal subgroup $S$ of $G$ generated
by all elements whose order is a power of $2$. Note that, since
$S_2$ has only one fixed point, $N_G(S_2)=G_{P}$ holds. If
$S=O(S)\rtimes S_2$, then $S_2$ is a Frobenius complement, and
hence it has a unique involutory element. Note that it is not
claimed that $S$ is a Frobenius group. In the exceptional cases,
if $C_G(S)$ is the centraliser of $S$ in $G$, then $G/C_G(S)$ is
an automorphism group of a group listed in Hering's theorem.

Therefore, the following result is obtained which together with
Lemma \ref{deuringshafarevic2} imply Theorem \ref{th1main}.

\begin{theorem}
\label{th1hering} Let $p=2$ and assume that {\rm (\ref{condstar})}
holds. Let $S$ be the subgroup $G$ of $\aut (\cX)$ of even
order$,$ generated by all elements whose order is a power of $2$.
One of the following holds$:$
\begin{enumerate}
\item[\rm(a)] $S$ isomorphic to one of the following groups$:$
\begin{equation}
\label{listahering} \PSL(2,n),\ \Sz(n),\,\PSU(3,n),\ {\mbox{\rm
SU}}(3,n),\,\,  \mbox{with $n =2^r\geq 4$};
\end{equation}
\item[\rm(b)] $S=O(S)\rtimes S_2$ with $S_2$  a  Sylow
$2$-subgroup of $G;$ here$,\ S_2$ is either a cyclic group or a
generalised quaternion group.
\item[\rm(c)] $S$ is a Sylow $2$-subgroup and $G=S\rtimes H,$ with
$H$ a subgroup of odd order.
\end{enumerate}
\end{theorem}

Using Theorem \ref{hering}, the action of $S$ on the set ${\Omega}$
of all points of $\cX$ fixed by some involution or, equivalently, on
the set consisting of all involutions in $G$ can be investigated.
Let $\bar{S}$ be the permutation group induced by $S$ on $\Omega$.
If $M$ is the kernel of the permutation representation of $S$, then
$|M|$ is odd, $\bar{S}=S/M$, and $\bar{S}_2=S_2M/M$ is a Sylow
$2$-subgroup of $\bar{S}$. Here, $\bar{S}$ is transitive on
$\Omega$; also, $\bar{S}_2$ fixes ${P}$ and acts on
$\Omega\setminus\{P\}$ as a semi-regular permutation group. Further,
the kernel of the permutation representation of $G$ on $\Omega$ is
$N=C_G(S)$. If $S$ is isomorphic to one of the groups in
(\ref{listahering}), then $\bar{S}$ is $2$-transitive on $\Omega$,
and $|\Omega|=n+1,\,n^2+1,n^3+1$ according as $\bar{S}\cong
\PSL(2,n),\,\Sz(n),\,\PSU(3,n)$.


\section{Upper bound on the orders of solvable subgroups of $\K$--automorphism groups of
{$\lowercase{p}$}--rank zero curves}

\begin{theorem}
\label{solvablelarge} Suppose $p>2$ and let $G$ be a solvable
subgroup of $\aut (\cX)$ satisfying condition {\rm (\ref{condstar})}.
Assume that $G$ fixes no point$,$ and that
\begin{equation}
\label{solvablelargecond} \mbox{$|G|$ is divisible by $p^2$.}
\end{equation}
If $\cX$ has genus $g\geq 2,$ then $|G|\leq 24 g(g-1)$.
\end{theorem}

\begin{proof}
Since $p$ is odd,
$\cX$ has more than one Weierstrass points. So, if $g=2$, then
$|G|\leq 48$ and hence Theorem \ref{solvablelarge} is true.

Let $N$ be a minimal normal subgroup of $G$. Since $G$ is solvable,
$N$ is an elementary abelian group of order $d^r$ with a prime $d$.
If $d=p$, then (\ref{condstar})  implies that $N$ fixes a unique
point $P$. But then $G$ itself must fix $P$, contradicting one of
the hypotheses.
So, $p\neq d$.
Theorem \ref{gow} yields that the Sylow $p$-subgroups of $G$ are
cyclic. As $d\neq p$, the Sylow $p$-subgroups of $G$ and $G/N$ are
isomorphic,  while (\ref{solvablelargecond}) remains valid for
$G/N$.

This suggests  that $\bar{G}=G/N$, viewed as a subgroup of
$\aut(\cY)$ with $\cY=\cX/G$, should be investigated. The aim is
to show that $\bar{G}$ also satisfies (\ref{condstar}).

The point ${\bar P}\in \cY$ lying below $P$ is fixed by a Sylow
$p$-subgroup of $\bar{G}$. Since $G/N$ can be viewed as a subgroup
$\aut(\cY)$, this and (\ref{solvablelargecond}) rule out the
possibility that $\cY$ is elliptic. Similarly, $\cY$ cannot be
rational. In fact, if $\cY$ were rational, then no $p$-subgroup of
$G/N$ and hence of $G$ would be  a cyclic group of order $p^i$
with $i\geq 2$,
contradicting Theorem \ref{gow}.

 Therefore, $\cY$  has genus $\bar{g}\geq 2$. Assume that
Theorem \ref{solvablelarge} does not hold, and choose a minimal
counterexample with genus $g$ as small as possible. Then
$|G|>24g(g-1)$, but Theorem \ref{solvablelarge} holds for all
genera $g'$ with $2\leq g'<g$. (\ref{eq1}) applied to $N$ gives
$2g-2\geq |N|(2\bar{g}-2),$ whence
$$
|G|>12g(2g-2)\ge |N|\,12g(2\bar{g}-2).
$$
Therefore $|G/N|>24\bar {g}(\bar{g}-1)$. By the divisibility
condition (\ref{solvablelargecond}) on $|G|$, the order of
$\bar{G}=G/N$ is not a prime. Hence $\bar{G}$ is a solvable group.
Since $d\neq p$, so $\bar{G}$ satisfies (\ref{solvablelargecond}).
On the other hand, $\bar{G}$ can be regarded as a subgroup of $\aut
(\cY)$.

To show that $\bar{G}$ also satisfies condition (\ref{condstar}),
let $\bar{\ga}$ be any element of $\bar{G}$ of order $p$. Choose
an element $\ga$ in $G$ whose image is $\bar{\ga}$ under the
natural homomorphism $G\mapsto G/N$. By (\ref{condstar}), $\cX$ has
a unique point $P$ fixed by $\ga$. Then, $\bar{\ga}$ fixes the
point ${\bar{P}}$ lying under $P$ in the covering $\cX\to \cY$.

Assume on the contrary that $\bar{\ga}$ fixes another point of
$\cY$, say $\bar{Q}$. Then the set of points of $\cY$ lying over
$\bar{Q}$ is preserved by $\ga$. The number of such points is
prime to $p$ since their number divides $|N|$. But then $\ga$ must
fix one of these points, contradicting (\ref{condstar}).

Since $G$ is a minimal counterexample, $\bar{G}$ fixes a point of
$\cY$. Equivalently, an orbit of $N$ is preserved by $G$. Such an
orbit consists of points of $\cX$ fixed by $p$-elements of $G$.
Since any Sylow $p$-subgroup $S_p$ of $G$ has only one fixed point, so $|S_p|\leq |N|-1$. This
together with (\ref{solvablelargecond}) implies that $|N|\geq 10$.
As $g-1\geq |N|(\bar{g}-1)$, it follows that $g\geq 11$.

Note that
$|\bar{G}|=|\bar{G}_{{\bar P}}|=|\bar{G}_{{\bar P}}^{(1)}|\,|\bar{H}|$, with
$\bar{H}$  a cyclic group of order prime to $p$. By Stichtenoth's
bound, $|\bar{H}|\leq4\bar{g}+2$. Let $\bar{S}_p$ be the Sylow
$p$-subgroup of $\bar{G}$; then $\bar{S}_p=\bar{G}_{{\bar P}}^{(1)}$. Since $p\nmid |N|$, so $\bar{S}_p$ is isomorphic to
$S_p$.

By Theorem \ref{gow}, $S_p$ is cyclic. By Nakajima's bound,
$|\bar{G}_{{\bar P}}^{(1)}|\leq 4\bar{g}+4$. Therefore,
$$
|\bar{G}|\leq (4\bar{g}+4)(4\bar{g}+2),
$$ whence $|G|\leq
(4\bar{g}+4)(4\bar{g}+2)|N|$. Since $g-1\geq |N|(\bar{g}-1)$,
$\bar{g}\geq 2$ and $|N|\geq 10$, it follows that
\begin{eqnarray*}
|G|& \leq & 16(\bar{g}+1)(\bar{g}+\ha)|N|\\
&\leq & \frac{16}{10}\cdot
\frac{\bar{g}+1}{\bar{g}-1}(\bar{g}-1)|N|\cdot
\frac{\bar{g}+\ha}{\bar{g}-1}(\bar{g}-1)|N|\\
&\leq &
 \frac{8}{5}\cdot\frac{2\bar{g}^2+3
\bar{g}+1}{2(\bar{g}-1)^2}(g-1)^2\\
&\leq & 12(g-1)^2<24g(g-1).
\end{eqnarray*}
But then $G$ is not a counterexample.
\end{proof}
Under the hypothesis that $16$ divides $|G|$, the above proof with
some modifications also works for $p=2$. Here an alternative proof
is given for $p=2$ which does not require that hypothesis.
\begin{theorem}
\label{th2hering}  Let $p=2,\ g\geq 2,$ and assume that
{\rm (\ref{condstar})}  holds. If  $G$ is a solvable subgroup of $\aut
(\cX)$ fixing no point of $\cX$ then $|G| \leq 24g^2$.
\end{theorem}

\begin{proof}
Since $84(g-1)<24g^2$, the group $G$ may be supposed to be one of
the exceptions in \cite{stichtenoth1973I}, Satz 3. Since $G$ has
no fixed point, from Theorem \ref{th1hering}, $S=O(S)\rtimes S_2$,
and $O(S)$ acts transitively on the set $\Omega$ consisting of all
points fixed by involutory elements in $G$. In particular,
$\Omega$ is the unique non-tame orbit of $G$.  Let $P\in \Omega$.
By assumption (\ref{condstar}), $G_{P}^{(1)}$ is a Sylow
$2$-subgroup $S_2$ of $G$. From $|O(S)|=|\Omega||O(S)_{P}|$, it follows
\begin{equation}
\label{eqth2hering}|G|=\frac{|G_{P}||O(S)|}{|O(S)_{P}|}.
\end{equation}
In particular, $\Omega$ has odd size. As $O(S)$ is transitive on
$\Omega$ and $N_G(S_2)=G_{P}$, each element in $G$ is the product
of an element fixing $P$ and an element from $O(S),$ that is,
$G=G_{P}O(S).$ Also, as $O(S)$ is a characteristic subgroup of $S$
and $S\trianglelefteq G$, so $O(S)\trianglelefteq G$. Hence
$\bar{G}=G/O(S)$ is a factor group which can be viewed as a
subgroup of
$\aut\, (\cY)$
where $\cY=\cX/O(S)$.

Assume that $\cY$ is rational. Then every $2$-subgroup of
$\bar{G}$ is elementary abelian. On the other hand, as $S_2\cong
S_2O(S)/O(S)$, from Theorem \ref{gow} the factor group
$S_2O(S)/O(S)$ has only one involutory element. Thus,
$|S_2O(S)/O(S)|=2$, whence $|S_2|=2$. This implies that
$G=O(G)\rtimes S_2$. Hence
$$
|G|=2|O(G)|\leq 168(g-1)\leq 24g^2
$$
for $g\geq 6$.

To show that this holds true for $2\leq g \leq 5$, some results
from \cite{stichtenoth1973I} are needed, see the proof of Theorem
3: If $G$ has more than three short orbits, then $|G|\leq
12(g-1)<24g^2$. Since $S_2$ has exactly one fixed point, it
preserves only one orbit of $O(G)$. Therefore,
either $O(G)$ (and hence $G$) has only one short orbit, or $G$ has $2$ tame orbits and $1$ wild orbit.
Since $p=2$, in the latter case $|G|\leq 84(g-1)$ as it follows from another
result shown in the proof of the above mentioned Theorem 3.
Therefore, $|G|< 24g^2$.

Assume that $\Omega$ is the unique short orbit of $G$. Then
$|\Omega|$ divides $2g-2$. Therefore,
$$
|G|\le |\Omega||G_P|\le (2g-2)|S_2||H|\le (4g-4)|H|,
$$
where $H$ is the complement of $S_2$ in $G_P$. Then Stichtenoth's bound
$|H|\le 4g+2$ yields $|G|<24g^2$.

Assume that $\cY$ is elliptic.
{}From (\ref{eq2}),
$$
2g-2=|O(S)|\sum_{i=1}^k\left(\frac{|O(S)_{P_i}|-1}{|O(S)_{P_i}|}\right),
$$
where $P_1,\ldots,P_k$ are representatives of the short orbits of $O(S)$.
Note that $k\ge 1$, as otherwise $g=1$.
Hence,
$|O(S)|<4(g-1)$. As $G_PO(S)/O(S)$ fixes a point of $\cY$, it follows that $|G|=
|G_PO(S)| \le 24 |O(S)| \le 96(g-1)\leq 24g^2$.

We are left with the case that the genus of the quotient curve
$\cY=\cX/O(S)$ is greater than $1$. From (\ref{eq2}), $|O(S)|\leq
g-1$. Hence, by (\ref{eqth2hering}),
$$
|G|\leq (g-1)|G_{P}|.
$$
Assume that the quotient curve $\cZ=\cX/{S_2}$ has genus $g'>1.$
Then $|S_2|(g'-1)<g-1$ by (\ref{eq1}). Since $G_{P}/S_2$ is a tame
subgroup of $\aut(\cZ)$ fixing the point of $\cZ$ lying under $P$,
from Stichtenoth's bound, $|G_{P}/S_2|\leq 4g'+2$. Hence,
$$
|G_{P}|<\frac{4g'+2}{g'-1}(g'-1)|S_2|\leq 10(g-1).
$$
Since $|\Omega|\leq |O(S)|\leq g-1$, this implies that
$|G|<10(g-1)^2<24g^2$.

If $\cZ$ is elliptic, then $|G_{P}/S_{2}|\leq 24$, while
$|S_{2}|\leq g$  by Stichtenoth's bound. Therefore, $|G|\leq
24g(g-1)<24g^2$.

So, take $\cZ$ to be rational. Let $R\in \Omega$ be a point
distinct from $P$. If the stabiliser $G_{P,R}$ of $R$ in $G_{P}$
is trivial, then $|G_{P}|<|\Omega|$, and hence
$$
|G|=|G_{P}||\Omega|< |\Omega|^2\leq |O(S)|^2\leq (g-1)^2.
$$
Therefore, let $|G_{P,R}|>1$ for every $R\in \Omega$.
We show that $G$ acts on $\Omega$ as a $2$-transitive permutation
group $\bar{G}$.

Let $\Omega_0=\{P\},\Omega_1,\ldots
\Omega_r$ with $r\geq 1$  denote the orbits of $S_2$
contained in $\Omega$. Then, $\Omega=\bigcup_{i=0}^r \Omega_i$. To
prove that $G$ acts $2$-transitively on $\Omega$, it suffices to
show that $r=1$. For any $i$ with $1\leq i \leq r$, take a point
$R\in \Omega_i$. By hypothesis, $R$ is fixed by an element $\ga\in
G_{P}$ whose order $m$ is a prime different from $p$. Since
$|G_P|=|S_2||H|$ and $m$ divides $|G_P|$, this implies that
$m$ must divide $|H|$. By the Sylow theorem, there is a subgroup
$H'$ conjugate to $H$ in $G_{P}$ which contains $\ga$; here, $\ga$
preserves $\Omega_i$. Since the quotient curve
$\cZ$ is rational, $\ga$ fixes at most two orbits
of $S_2$. Therefore, $\Omega_0$ and $\Omega_i$ are the
orbits preserved by $\ga$. As $H'$ is abelian and $\ga\in H'$,
this yields that $H'$ either preserves both $\Omega_0$ and
$\Omega_i$ or interchanges them. The latter case cannot actually
occur as $H'$ preserves $\Omega_0$. So, the orbits $\Omega_0$ and
$\Omega_i$ are also the only orbits of $S_2$ which are
fixed by $H'$. Since $G_{P}=S_2\rtimes H'$, this implies
that the whole group $G_{P}$ fixes $\Omega_i$. As $i$ can be any
integer between $1$ and $r$, it follows that $G_{P}$ fixes each of
the orbits $\Omega_0,\Omega_1,\ldots,\Omega_r$. Hence, either
$r=1$ or $G_{P}$ preserves at least three orbits of $S_2$. The latter case cannot actually occur, as the quotient
curve $\cZ$  is rational. Therefore $r=1$. Also,
the size of $\Omega$ is of the form $q+1$ with $q=S_2$;
in particular, $q$ is a power of $2$, and $\Omega$ comprises $P$
together with all points in the unique non-trivial orbit of $S_2$.

The latter assertion implies that $|\Omega|=|S_2|+1=q+1$, $q=2^r$.
Let ${\bar G}$ denote the $2$-transitive permutation group induced by $G$ on $\Omega$.
As $G$ is solvable, ${\bar G}$ admits
an elementary abelian $w$-group acting on $\Omega$ as a sharply transitive permutation group.
Hence, $|\Omega|=w^k$ with an odd prime $w$. {}From $q+1=w^k$,
either $k=1$, or $k=2$, $q=8$, $w=3$,
see \cite{mihailescu2004}.
Let $M$ be the kernel of the permutation representation of $G$ on
$\Omega$.

It may be that $M$ is trivial, that is $G=\bar{G}$, and this
possibility is investigated first. Assume that $k=1$. Then
$G=\AGL(1,w)$ and $G$ is sharply $2$-transitive on $\Omega$. Let
$N$ be the normal subgroup of $G$ of order $w$. If the quotient
curve $\cU=\cX/N$ has genus $g'\ge 2$, then
$$
2g-2\geq w(2g'-2)\ge 2w,
$$
and hence
$$
|G|=|\Omega|(|\Omega|-1)=w(w-1)\leq (g-1)(g-2)<24g^2.
$$

If $\cU$ is elliptic, then $|S_2|$ is either $2$, or $4$, or $8$.
The last case cannot occur because $9$ is not a prime number. If
$|S_2|=2$, then $w=3$, and hence $|G|=6$, while $|S_2|=4$ occurs
when $w=5$ and hence $|G|=20$.

If $\cU$ is rational, then $S_2$, viewed as a subgroup of
$\aut(\cU)$ must be an elementary abelian group. On the other
hand, $S_2$ is the $1$-point stabiliser of $\AGL(1,w)$ which is
cyclic. It follows that $|S_2|=2$. Therefore, $|o|=3$, whence
$|G|=6$.

If $k=2,w=3,q=8$, then $G$ has an elementary abelian normal subgroup $N$ of order $9$.
Since $G$ is a subgroup of the
normalizer of $N$, it follows that $G_P$ is a subgroup of
$\mathrm{GL}(2,3)$ whose order is divisible by eight but not by
sixteen. Since $|\mathrm{GL}(2,3)|=48$, either $|G_P|=8$, or
$|G_P|=24$. In the former case, $|G|=72$ and hence $|G|<24g^2$. In
the latter case, $|G|=216,$ which is greater than $24g^2$ only for
$g=2$. But this does not actually occur,
as $\cX$ has more than one Weierstrass -point and hence
the assertion holds for $g=2$.

Suppose $M$ is non-trivial. Then $M$ is cyclic of odd order. Let
$\cU$ be the quotient curve $\cX/M$. Then $\bar{G}=G/M$ can be
viewed as a subgroup of $\aut(\cU)$. The Sylow $2$-subgroup
$\bar{S}_2=S_2M/M$ of $\bar{G}$ is isomorphic to $S_2$.
Arguing as in the proof of Theorem \ref{solvablelarge} shows that
$\bar{G}$ has property (\ref{condstar}).

{}From above, if the genus $\bar{g}$ of $\cU$ is at least $2$,
then $|\bar{G}|\leq 24 \bar{g}^2$; in consequence, $|G|\leq
24|M|\bar{g}^2$. Since $|\Omega|\geq 3$, from (\ref{eq2}),
$$2g-2\geq |M|(2\bar{g}-2)+3(|M|-1)\geq 2|M|\bar{g}-2.$$
This shows that $|G|\leq 24 g \bar{g}<24g^2$.

If $\cU$ is elliptic, then $G/M$ has order at most $24$ and the
assertion holds for $|M|=3$. Furthermore, ${S_2}\cong S_2M/M$
consists of at most $8$ elements. Assume that $|M|>3$. Then
$|S_2|\leq 2(|M|-1)$. From $2g-2\geq (|M|-1)|\Omega|$,
$$
\mid G\mid =|G_P||\Omega|\le
|S_2|(4g+2)|\Omega|=2(4g+2)(|M|-1)|\Omega|\le 4(g-1)(4g+2)
$$
which is smaller than $24g^2$.

Finally, if $\cU$ is rational, then $|S_2|=2$, and the argument above gives
$$
|G| \le | S_2|(4g+2)(|S_2 |+1) \le 6(4g+2)<24g^2.
$$
%
\end{proof}
\section{Large simple $\K$--automorphism groups of {$\lowercase{p}$}--rank zero curves}
Throughout this and the next sections, the hypotheses of Theorem
\ref{th1main} are assumed, and notation as used in Sections
\ref{sec1}, \ref{sec2} and \ref{sectionstructure} is maintained.
Our purpose is to prove the following result.

\begin{theorem}  
\label{th2main}
The possibilities in Theorem \ref{th1main} are the
following.

\begin{enumerate}
        \item[\rm(a)]\quad $S$ fixes a point of $\cX$.
        \item[\rm(b)]\quad $S$ is solvable and $|S|\leq 24g^2$.
        \item[\rm(c)]\quad $S\cong \PSL(2,n)$ and if both $|S|>6(g-1)$ and {\rm (\ref{eqfinal})} hold, then
        $g=\ha(t-1)(n-1)$ with $t|(n+1).$
        \item[\rm(d)]\quad  $S\cong \PSU(3,n)$ and
if both $|S|> 6(g-1)$ and {\rm (\ref{eqfinal})} hold, then either
        $$g=\ha\,(n-1)(t(n+1)^2-(n^2+n+1))$$ with
    $t|(n^2-n+1)/\mu$, or
$$
g=\ha(n-1)\left(\frac{t(n^3+1)}{\mu}-(n^2+n+1)\right),
$$
with $t\mid (n+1)$.    In the former case, if $t=1$ then $\cX$ is
as in {\em (II)}.
        \item[\rm(e)]\quad $S\cong
{\rm{SU}}(3,n)$ with $3|(n+1)$ and if both $|S|> 6(g-1)$ and {\rm
(\ref{eqfinal})} hold, then either
$$g=\ha\,(n-1)[3t(n+1)^2-(n^2+n+1)]$$  with $t|(n^2-n+1)/3$, or
$$
g=\ha (n-1)\left({t(n^3+1)}-(n^2+n+1)\right),
$$
with $t\mid (n+1)$.
        \item[\rm(f)]\quad $S\cong \Sz(n)$ and
if both $|S|> 6(g-1)$ and {\rm (\ref{eqfinal})} hold, then
    either $$g\,=\ha\,[(t-1)(n^2-1)-2tn_0(n-1)]$$ with $t|(n
    +2n_0+1)$, or $$g\,=\ha\,[(t-1)(n^2-1)+2tn_0(n-1)]$$ with $t|(n
    -2n_0+1)$; if $t=1$ in the latter case then $\cX$ is as in {\em (III)}.
\end{enumerate}
\end{theorem}
For the purpose of the proof, $S$ may be assumed to be one of
the four groups in (\ref{lingp}). From the discussion at the end
of Section \ref{sectionstructure}, $S$ has exactly one non-tame
orbit, namely $\Omega$, on which $S$ acts as one of the groups
$\PSL(2,q)$, $\PSU(3,n)$ and $\Sz(n)$ in its natural
$2$-transitive permutation representation. In particular, every
complement $H$ of $S_P^{(1)}$ in $S_P$ has just one more fixed
point in $\Omega$.

Before investigating such cases separately, we prove a few lemmas.
\begin{lemma}
\label{extra1} Let $P\in \Omega$. Assume that $|S|> 6(g-1)$. If
$\cX_1=\cX/S_P^{(1)}$ is rational, then $S$ has exactly two short
orbits, namely $\Omega$ and a tame orbit $\Delta$.
\end{lemma}
\begin{proof} By an argument going back to Hurwitz and adapted for any
groundfield by Stichtenoth \cite{stichtenoth1973I}, the hypothesis
$|S|> 6(g-1)$ implies that $S$ has at most one tame orbit.
On the other hand, assume that $\Omega$ is the unique short orbit
of $S$. Since $\cX_1$ is rational, from (\ref{eq1}) and
(\ref{condstar}),
\begin{equation}\label{formuladP1}
2g-2+2|S_P^{(1)}|=2|S_P^{(1)}|-2+|S_P^{(2)}|-1+\ldots=d_P-|S_P|+|S_P^{(1)}|.
\end{equation}
Therefore,
\begin{equation}\label{formuladP}
d_P=2g-2+|S_P^{(1)}|+|S_P|.
\end{equation}
Also, from (\ref{eq1}),
\begin{equation}
\label{eqrefinecaseiii} 2g-2=-2|S|+\deg\, D(\cX/S)= |\Omega|\,(
d_P
-2|S_{P}|).
\end{equation}
Hence $|\Omega|(|S_P|-|S_P^{(1)}|)=(|\Omega|-1)(2g-2).$ Since $S$
acts on $\Omega$ as one of the groups $\PSL(2,q)$, $\PSU(3,n)$ and
$\Sz(n)$ in its natural $2$-transitive permutation representation,
$|\gO|=|S_P^{(1)}|+1$ holds. Thus, $|\Omega|(|H|-1)=2g-2$ for a
complement $H$ of $S_P^{(1)}$ in $S_P$. {}From this
$|\Omega|^2|H|^2<16g^2$. Therefore,
$$|S|=|S_P^{(1)}||H||\Omega|=(|\Omega|-1)|H||\Omega|<16g^2,$$
a contradiction. \end{proof}
\begin{lemma}
\label{ex2} If {\rm (\ref{eqfinal})} holds then
 $\cX_1=\cX/S_P^{(1)}$ is rational.
\end{lemma}
\begin{proof}
Assume on the contrary that the genus $g'$ of $\cX_1$ is positive.
{}From (\ref{eq1}) applied to $S_P^{(1)}$,
\begin{equation}
\label{boundGPeq1} g\geq g'\,|S_P^{(1)}|.
\end{equation}
Take a complement $D$ of $S_P^{(1)}$ in $\aut(\cX)_P$. $D$
acts faithfully on $\cX_1$ as a subgroup of $\aut(\cX_1)$.
Furthermore, $D$ has at least two fixed points,
as $S$ is one of the groups in (\ref{lingp}),
namely the points ${\bar P}$ and ${\bar Q}$ lying under $P$ and a
point $Q\in \Omega\setminus\{P\}$  in the covering $\cX\mapsto
\cX_1$. Let $g''$ be the genus of the quotient curve $\cZ=\cX_1/D$
of $\cX_1$ with respect to $D$. From (\ref{eq2}) applied to $D$,
$$2g'-2\geq |D|(2g''-2)+2(|D|-1).$$ If $g''\geq 1$, this yields
$g'\geq |D|$ whence $g\geq |\aut(\cX)_P|$ by (\ref{boundGPeq1}).
Assume that $g''=0$. {}From (\ref{eq2}),
$$2g'=\sum_{i=1}^k(|D|-\ell_i)$$ with $\ell_1,\ldots,\ell_k$ being
the sizes of the short orbits of $D$ on $\cX_1$, other than the two trivial
ones $\{{\bar P}\}$ and $\{{\bar Q}\}$. Since $|D|$ is odd, because $p=2$,
each such short orbit has length at most $\thi |D|$. From this,
$g'\geq \thi |D|$ whence $g\geq \thi |\aut(\cX)_P|$  by
(\ref{boundGPeq1}). But this contradicts (\ref{eqfinal}).
\end{proof}

\begin{lemma}
\label{GQfpf} Assume that $|S|> 6(g-1)$ and that
$\cX_1=\cX/S_P^{(1)}$ is rational. Let $Q\in \Delta$. Then the
subgroups $S_{P}$ and $S_{Q}$ have trivial intersection, and $S_Q$
is a cyclic group whose order divides $q+1$. Also,
\begin{equation}
\label{Sticht13.11tris}
2g-2=\frac{|S|\,(|S_{P}|-|S_{P}^{(1)}|\,|S_{Q}|)} {|S_{Q}|(|S| -
|S_{P}|)}
\end{equation}
\end{lemma}
\begin{proof} Let $\ga\in S_{P}\cap S_{Q}$ be non-trivial. Then
$p\nmid {\rm{ord}}\, \ga$, and hence $\ga$ is in a complement $H$
of $S_P^{(1)}$. Therefore, $\ga$ fixes not only $P$ but another
point in $\Omega$, say $R$. Since $Q\not\in \Omega$, this shows
that $\ga$ has at least three fixed points. These points are in
three different orbits of $S_{P}^{(1)}$. Since the quotient curve
$\cX_1=\cX/S_P^{(1)}$ is rational, this implies that $\ga$ is
trivial, a contradiction. Hence $|S_{P}\cap S_{Q}|=1$. Therefore,
no non-trivial element of $S_Q$ fixes a point in $\Omega$. Since
$|\Omega|=q+1$, the second assertion follows.
Let
\begin{equation}
\label{eta} \eta=|S_Q|(d_P-|S_P|)-|S_P|.
\end{equation}
Then
\begin{equation}
\label{Sticht13.7} |S|=2(g-1)\frac{|S_P^{(1)}||H||S_Q|}{\eta},
\end{equation}
where $H$ is a complement of $S_P^{(1)}$ in $S_P$. Substituting
$d_P$ from (\ref{formuladP}) into (\ref{eta}) and then this into
$(\ref{Sticht13.7})$ gives (\ref{Sticht13.11tris}).
\end{proof}

    \subsection{Case $S\cong \PSL(2,q)$}
In this case $n=q$ and
$$|S|=q^3-q,\,\,|S_P|=q(q-1),\,\\\,|S_P^{(1)}|=q.$$
By Lemma \ref{ex2}, if (\ref{eqfinal}) holds then $\cX_1$ is
rational. Lemma \ref{GQfpf} implies that $|S_Q|$
divides $q+1$. Let $|S_Q|=(q+1)/t$ with $t|(q+1)$. Equation
(\ref{Sticht13.11tris}) gives
$2g=(t-1)(q-1)$.

    \subsection{Case $S\cong\PSU(3,n)$}
In this case, $q=n^3$ and
$$
|S|=(n^3+1)n^3(n^2-1)/\mu,\ |S_{P}|=n^3(n^2-1)/\mu,\
|S_{P}^{(1)}|=n^3,
$$
Also, by Lemma \ref{ex2}, (\ref{eqfinal}) implies that $\cX_1$ is
rational.
 By Lemmas \ref{backgr} and \ref{GQfpf}, two
cases are to be discussed, according as $|S_Q|$ divides either
$(n^2-n+1)/\mu$, or $n+1$. By (\ref{Sticht13.11tris}),
$$
2g=\frac{(n^3+1)(n^2-1)} {\mu |S_{Q}|}-(n^3-1).
$$

If $|S_Q|=(n^2-n+1)/(t\mu)$, then
\begin{equation}
\label{psubis} 2g=(n-1)(t(n+1)^2-(n^2+n+1)).
\end{equation}
If $t=1$,  then
$g=\ha n(n-1)$. Therefore, $t=1$ only occurs when (II) holds.

If $|S_Q|=(n+1)/t$, then
$$2g=(n-1)\left[\frac{(n^3+1)t}{\mu}-(n^2+n+1)\right]. $$

\subsection{Case $S\cong\SU(3,n)$ with $3|(n+1)$}
Let $\cY=\cX/Z(S)$ be the quotient curve of $\cX$ with respect to
the centre $Z(S)$ of $S$. By Lemma \ref{GQfpf}, the points in
$\Omega$ are the fixed points of the non-trivial elements in
$Z(S)$. The group $\bar{S}=S/Z(S)\cong \PSU(3,n)$ can be viewed as
a subgroup of $\aut(\cY)$. Since $|Z(S)|=3$, from (\ref{eq2})
$$2g-2=3(2\bar{g}-2)+2(n^3+1)$$ where $\bar{g}$ is the genus of
$\cY$. From this, if (\ref{eqfinal}) holds, then
$|\aut(\cY)_{\bar{P}}|$ is even and bigger than $3\bar{g}$ for the
point ${\bar P}$ lying under $P$ in the covering $\cX\to \cY$.
Also, $|{\bar S}|>6({\bar g}-1)$ holds. Hence, the preceding case
applies to $\cY$ whence the assertion follows for $\cX$.

\subsection{Case $S=\Sz(n)$}
In this case, $n=2n_0^2,\,n_0=2^r,\,r\geq 1$, and
$$
|S|=(n^2+1)n^2(n-1),\ |S_{P}|=n^2(n-1),\ |S_{P}^{(1)}|=n^2.
$$
Also, by Lemma \ref{ex2}, if (\ref{eqfinal}) holds then $\cX_1$ is
rational. From Lemmas \ref{backgr} and \ref{GQfpf}, there is an
odd integer $t$ such that either (A) or (B) holds, where
\begin{eqnarray*}
\mbox{(A)}\hspace*{8mm} & & 2g\,=\,(t-1)(n^2-1)-2tn_0(n-1), \quad
|S_{Q}|=(n +2n_0+1)/t;\\
\mbox{(B)}\hspace*{8mm} & & 2g\,=\, (t-1)(n^2-1)+2tn_0(n-1),\quad
|S_{Q}|=(n-2n_0+1)/t.
\end{eqnarray*}
In case (B), if $t=1$ then (III) holds.

\section{Non-solvable $\K$--automorphism groups of {$\lowercase{p}$}--rank zero curves}
\label{nonsolvableaut}

\begin{theorem}
\label{thmain} Let $\cX$ be a zero $2$-rank algebraic curve of
genus $g\geq 2$ defined over an algebraically closed groundfield
$\K$ of characteristic $2$.  If $\aut(\cX)$ is non-solvable then
its commutator group $\aut(\cX)'$ coincides with $S$ and hence is
one of the groups in {\rm (\ref{lingp})}.
\end{theorem}
\begin{proof} Since groups of odd order, as well as subgroups of $\aut(\cX)$ fixing a point of $\cX$, are
solvable, case (a) in Theorem \ref{th1main} holds. In particular,
$S$ is a perfect group, that is,  $S$ coincides with its commutator
group.  Then $S$ is a subgroup of the commutator subgroup
$\aut(\cX)'$ of $\aut(\cX)$.

On the other hand, since $S$ is a normal subgroup of $\aut(\cX)$,
the factor group $\aut(\cX)/S$ can be viewed as a $\K$-automorphism
group of the quotient curve $\cZ=\cX/S$. For a point $P$ in
$\Omega$, let ${\bar P}$ be the point of $\cZ$ lying under $P$.
Since $\Omega$ an $S$-orbit on $\cX$, the set of points lying over
${\bar P}$ coincides with $\Omega$. As $\Omega$ is also an
$\aut(\cX)$-orbit, this implies that the factor group $\aut(\cX)/S$
fixes ${\bar P}$. Furthermore, $\aut(\cX)/S$ is tame as
$\aut(\cX)/S$ has odd order. Hence, $\aut(\cX)/S$ is cyclic.
Therefore, $S$ contains $\aut(\cX)'$.
\end{proof}

For a non-solvable group $\aut(\cX)$, the following theorem
describes the structure of $\aut(\cX)$ in terms of $\aut(\cX)'$ and
the kernel $N$ of the permutation representation of $\aut(\cX)$ on
$\Omega$.
\begin{theorem}
\label{thmaincor1} Let $\cX$ be a zero $2$-rank algebraic curve of
genus $g\geq 2$ defined over an algebraically closed groundfield
$\K$ of characteristic $2$. If $\aut(\cX)$ is non-solvable then
$\aut(\cX)'=S$ and one of the following cases occur for a cyclic group of odd order $N$:
\begin{itemize}
    \item[\rm(i)] $\aut(\cX)'=\PSL(2,n)$ and $\aut(\cX)=\PSL(2,n)\times
    N$;
    \item[\rm(ii)] $\aut(\cX)'=\Sz(n)$ and $\aut(\cX)=\Sz(n)\times
    N$;
    \item[\rm(iii)] $\aut(\cX)'=\PSU(3,n)$ and either
    \begin{itemize}
    \item[\rm(iii)(a)] $\aut(\cX)/N=\PSU(3,n)$ and
    $\aut(\cX)=\PSU(3,n)\times N$,
    or
    \item[\rm(iii)(b)] $\aut(\cX)/N=\PGU(3,n)$ and $\PSU(3,n)\times N$
    is a subgroup of index $3$ of $\aut(\cX)$.
    \end{itemize}
    \item[\rm(iv)] $\aut(\cX)'=\SU(3,n)$ and either

    \begin{itemize}
    \item[\rm(iv)(a)] $\aut(\cX)/N=\PSU(3,n)$ and
    $\aut(\cX)=\SU(3,n)N$,
    or
    \item[\rm(iv)(b)] $\aut(\cX)/N=\PGU(3,n)$ and $\SU(3,n)N$
    is a subgroup of index $3$ of $\aut(\cX)$;
    \end{itemize}
\end{itemize}
where $n\geq 4$ is a power of $2$.
\end{theorem}
\begin{proof}
Consider the quotient curve $\cZ=\cX/N$ where $N$ is the subgroup of
$\aut(\cX)$ fixing $\Omega$ pointwise. Then $\aut(\cX)/N$ can be
viewed as an automorphism group of the quotient curve $\cX/N$. From
Lemma \ref{stabi} applied to $L=SN/N$, it follows that
$$
\aut(\cX)/N\cong \left\{\begin{array}{l}\PSL(2,n),\\
\PSU(3,n),\\ \PGU(3,n),\quad {\rm gcd}(3,n+1)=3,\\
\Sz(n).\end{array}\right.
$$
Furthermore,
$$
SN/N\cong \left\{\begin{array}{l}\PSL(2,n),\\
\PSU(3,n),\\ \Sz(n).\end{array}\right.
$$
{}From the third isomorphism theorem,
$$
\aut(\cX)/SN\cong (\aut(X)/N)/(SN/N).
$$
This implies that either $\aut(\cX)=SN$ or $\aut(\cX)/N\cong
\PGU(3,n)$ and $SN/N\cong \PSU(3,n)$  with ${\rm gcd}(3,n+1)=3$.
Note that  $S$ and $N$ have non-trivial intersection only for
$S=\SU(3,n)$ when $|S\cap N|=3$. Hence, if $\aut(\cX)=SN$ holds,
then one of the cases (i), (ii), (iiia), (iva) occurs. Furthermore,
if $\aut(\cX)/N\cong \PGU(3,n)$ then $SN$ has index $3$ in
$\aut(\cX)$, and either (iiib) or (ivb) holds.
\end{proof}

Finally, Theorem \ref{th1maincor} is a consequence of Theorems
\ref{th2main}, \ref{thmain} and \ref{thmaincor1} together with the
following lemma.
\begin{lemma}
\label{global} Assume that $|\aut(\cX)|>24g^2$. If $\aut(\cX)$ is
non-solvable, then both conditions {\rm (\ref{eqfinal})} and {\rm
(\ref{eqfinalbis})} are satisfied.
\end{lemma}
\begin{proof} Take a point $P\in \Omega$, and assume on the contrary that
 (\ref{eqfinal}) does not hold.
Since $|\aut(\cX)_P|$ is even, this only occurs when
$|\aut(\cX)_P|\leq 3g$. Then clearly $|S_P^{(1)}|\leq 3g$ holds.
This and $|\Omega|=|S_P^{(1)}|+1$ imply that $|\Omega|\leq
(3g+1)$. Therefore,
$$|\aut(\cX)|=|\aut(\cX)_P|\,|\Omega|\leq 3g(3g+1),$$ a contradiction
with $\aut(\cX)|\geq 24g^2$.

To prove that  (\ref{eqfinalbis}) holds, note that by Theorem
\ref{thmaincor1}
\begin{equation}\label{last}
|S|\ge \frac{|\aut(\cX)|}{3|N|}> \frac{8g^2}{|N|},
\end{equation}
 with $N$ a cyclic group of odd order
fixing $|\Omega|=|S_P^{(1)}|+1\ge 3$ points of $\cX$. Then, by
(\ref{eq2}),
$$
2g-2\ge -2|N|+|\Omega|(|N|-1)\ge (|N|-1)(|\Omega|-2)-2.
$$
Since $|N|\geq 3$ and $|\Omega|\geq 5$, this yields that $2g\geq
\tth |N|(|\Omega|-2)\geq 2|N|$. From this and (\ref{last}), we
obtain that $|S|>8g>6(g-1)$.
\end{proof}

\section{Some examples}
As pointed out in Introduction, each case in Theorem \ref{th1main}
(a) occurs, examples being the four curves (I)-(IV). To complete
the discussion after Theorem \ref{th1main}, we exhibit further
examples for (i) and for the first case in both (ii) and (iv). As
far as we know, no more examples exist.

In case (i), Stichtenoth's result, see \cite{stichtenoth1973II},
Satz 7, determines all examples of curves $\cX$ with an absolutely
irreducible plane model $\cC$
\begin{equation}
\label{a(y)=b(x)} \bv(A(Y)+B(X)),
\end{equation}
where
\begin{enumerate}
\item[(I)] $\deg\, \cC\geq 4$;

\item[(II)] $ A(Y)=a_n Y^{2^n}+a_{n-1} Y^{2^{n-1}} +\ldots+ a_0
Y,\quad a_j\in K,\, a_0,\ a_n\neq 0;$
\item[(III)]
 $B(X)=b_mX^m+b_{m-1}X^{m-1}+\ldots+ b_1X+b_0, \quad
b_j\in K,\ b_m\neq 0;$

\item[(IV)] $ m$ is odd;

\item[(V)] $ n\geq 1, \ m\geq 3.$
\end{enumerate}
Such a curve $\cX$ has genus $g=\ha\,(m-1)(p^n-1),$ and
Stichtenoth's result for $p=2$ is as follows.
\begin{proposition}
\label{StichtarchIISatz7} If $\cC$ is a curve of type
{\rm{(\ref{a(y)=b(x)})}}$,$ then $\Aut(\cX)$ fixes a point, except
in the following two cases:
\begin{enumerate}
\item[\rm(i)]
\begin{enumerate}
    \item[\rm(a)] $\cC=\bv(Y^{2^n}+Y+X^m),$ with $m<2^n,\, 2^n\equiv-1
\pmod{m};$
    \item[\rm(b)] $\Aut(\cX)$ contains a cyclic subgroup
$C_m$ of order $m$ such that $\Aut(\cX)\cong \PSL(2,2^n)\times
C_m$.
\end{enumerate}
\item[\rm(ii)]
\begin{enumerate}
\item[\rm(a)]
 $\cX=\cC= \bv(Y^{2^n}+Y +X^{2^n+1}),$ the Hermitian curve;
 \item[\rm(b)] $\Aut(\cX)\cong \PGU(3,2^n).$
\end{enumerate}
\end{enumerate}
\end{proposition}


\begin{proposition}
\label{suex} Let $3|(n+1)$. For every divisor $t$ of
$(n^2-n+1)/3$, the non-singular model of the plane curve
\begin{equation}
\label{suexeq} \bv(Y^{(n^2-n+1)/t}-X^{n^3}-X+(X^n+X)^{n^2-n+1})
\end{equation}
of genus $g=\ha\,(n-1)(3t(n+1)^2-(n^2+n+1))$ has a
$\K$-automorphism group isomorphic to $\SU(3,n)$. Furthermore, the
non-singular model of the plane curve
\begin{equation}
\label{suexeqbis}
\bv(Y^{(n^2-n+1)/(3t)}-X^{n^3}-X+(X^n+X)^{n^2-n+1})
\end{equation}
of genus $g=\ha\,(n-1)(t(n+1)^2-(n^2+n+1))$ has a
$\K$-automorphism group isomorphic to $\PSU(3,n)$.
\end{proposition}
\begin{proof}
To show the first part, let $\cY$ be a non-singular model of the
curve (\ref{suexeq}). Note that $\cY$ can be regarded as the
quotient curve of the curve $\cX$ given in (IV) with respect to
the subgroup $H$ of $\aut(\cX)$ generated by the automorphism of
equation $x'=x,\,y'=\gl y$ where $\gl\in \K$ is a primitive
$t$-root of unity. The fixed points of the non-trivial
automorphisms in $H$ are the $n^3+1$ points of $\cX$ arising from
the $n^3+1$ branches of the plane curve (\ref{suexeq}) centred at
points on the $X$-axis. Since $t$ is odd and $\cX$ has genus
$\ha\,(n^3+1)(n^2-2)$, from (\ref{eq2})
$g=\ha\,(n-1)(3t(n+1)^2-(n^2+n+1))$ follows. The normalizer of $H$
in $\aut(\cX)$ contains a subgroup isomorphic to $\SU(3,n)$ with
intersects $H$ trivially. This completes the proof of the first
part.

To show the second part, let $\cZ$ be a a non-singular model of
the curve (\ref{suexeqbis}). This time, $\cZ$ is regarded as the
quotient curve of the curve $\cY$ with respect to the subgroup $H$
of $\aut(\cY)$ generated by the automorphism of equation
$x'=x,\,y'=\gep y$ where $\gep\in \K$ is a primitive third root of
unity. As before, the fixed points of the non-trivial
automorphisms in $H$ are the $n^3+1$ points of $\cY$ arising from
the $n^3+1$ branches of the plane curve (\ref{suexeq}) centred at
points on the $X$-axis. From (\ref{eq2}) applied to $\cY$ and $G$,
we get $g=\ha\,(n-1)(t(n+1)^2-(n^2+n+1))$. The centralizer of $H$
in $\aut(\cY)$ is a subgroup $T$ isomorphic to $\SU(3,n)$. Since
$\SU(3,n)/Z(\SU(3,n))\cong \PSU(3,n)$, this shows that $\aut(\cZ)$
has a subgroup isomorphic to $\PSU(3,n)$.
\end{proof}
\begin{proposition}
\label{psuex} Let $3|(n-1)$. For every divisor $t$ of $n^2-n+1$,
the non-singular model $\cX$ of the plane curve
\begin{equation}
\label{psuexeq} \bv(Y^{(n^2-n+1)/t}-X^{n^3}-X+(X^n+X)^{n^2-n+1})
\end{equation}
of genus $g=\ha\,(n-1)(3t(n+1)^2-(n^2+n+1))$ has a
$\K$-automorphism group isomorphic to $\PSU(3,n)$.
\end{proposition}
\begin{proof} For $3|(n-1)$, $\SU(3,n)=\PSU(3,n)$ holds, and the proof can be
carried out following the first part in the proof of Proposition
\ref{suex}.
\end{proof}
Our final remark is that the curves discussed in the last two
propositions are investigated in \cite{segovia}.

\vspace{0,5cm}\noindent {\em Authors' addresses}:

\vspace{0.2 cm} \noindent Massimo GIULIETTI \\
Dipartimento di Matematica e Informatica
\\ Universit\`a degli Studi di Perugia \\ Via Vanvitelli, 1 \\
06123 Perugia
(Italy).\\
 E--mail: {\tt giuliet@dipmat.unipg.it}

\vspace{0.2cm}\noindent G\'abor KORCHM\'AROS\\ Dipartimento di
Matematica\\ Universit\`a della Basilicata\\ Contrada Macchia
Romana\\ 85100 Potenza (Italy).\\E--mail: {\tt
gabor.korchmaros@unibas.it }


\begin{thebibliography}{99}



\bibitem{cakcak-ozbudak2004}
E.~\c{C}ak\c{c}ak and F.~\"Ozbudak,  Subfields of the function field
of the Deligne--Lusztig curve of Ree type. {Acta Arith.} {\bf
115} (2004), 133--180.



\bibitem{garcia-stichtenoth-xing2001} A.~Garcia, H.~Stichtenoth and
C.P.~Xing, On subfields of the Hermitian function field.
{Compositio Math.} {\bf  120} (2000),  137--170.



\bibitem{giulietti-korchmaros2007} M.~Giulietti and G.~Korchm\'aros,
On large automorphism groups of algebraic curves in positive
characteristic. Available at http://arxiv.org/abs/0706.2320.

\bibitem{segovia} M.~Giulietti and G.~Korchm\'aros,
A new family of maximal curves over a finite field. Available at
http://arxiv.org/abs/0711.0445.


\bibitem{giulietti-korchmaros-torres2004} M.~Giulietti,
G.~Korchm\'aros and F.~Torres,  Quotient curves of the
Deligne--Lusztig curve of Suzuki type. {Acta Arith.} {\bf 122}
(2006), 245--274.

\bibitem{greenberg1974} L.~Greenberg, Maximal groups and
signatures. { Ann. of Math. Studies} {\bf{79}} (1974), 207-226.

\bibitem{gow1975} R.~Gow,  Some $TI$ subgroups of solvable groups. {J. London
Math. Soc.} {\bf 12} (1976), 285--286.



\bibitem{hansen1993} J.P. Hansen and J.P. Pedersen, Automorphism
group of Ree type, Deligne-Lusztig curves and function fields. { J.
Reine Angew. Math.} {\bf 440} (1993), 99--109.



\bibitem{henn1978} H.W. Henn,  Funktionenk\"orper mit gro$\beta$er
Automorphismengruppe. { J. Reine Angew. Math.} {\bf 302} (1978),
96--115.

\bibitem{hering1972} C.~Hering,  On subgroups with trivial normalizer
intersection. {J. Algebra} {\bf 20} (1972), 622--629.

\bibitem{hirschfeld1985} J.W.P.~Hirschfeld,
\emph{Finite Projective Spaces of Three Dimensions}, Oxford Univ.
Press, Oxford, 1985, x+316 pp.

\bibitem{hirschfeld1998} J.W.P. Hirschfeld, {\em Projective Geometries Over Finite
Fields}, Second edition, Oxford University Press, Oxford, 1998, xiv+555 pp.

\bibitem{hirschfeld-korchmaros-torres2008} J.W.P.~Hirschfeld,
G.~Korchm\'aros and F.~Torres \emph{Algebraic Curves Over a Finite
Field}, Princeton Univ. Press, Princeton and Oxford, 2008, xx+696
pp.

\bibitem{huppertI1967} B.~Huppert, \emph{Endliche Gruppen. I}, Grundlehren der Mathematischen
Wissen\-schaften {\bf 134}, Springer, Berlin, 1967, xii+793 pp.


\bibitem{huppertblackburn3-1982}
B.~Huppert and B.N.~Blackburn, \emph{Finite groups. III},
Grundlehren der Mathematischen Wissenschaften {\bf 243}, Springer,
Berlin, 1982, ix+454 pp.


\bibitem{iwasawatamagawa1951} K.~Iwasawa and T.~Tamagawa,  On the group of
automorphisms of a function field. {J. Math. Soc. Japan} {\bf
3} (1951), 137--147.


\bibitem{iwasawatamagawa1951c} K.~Iwasawa and T.~Tamagawa, Correction:
On the group of automorphisms of a function field. {J. Math.
Soc. Japan} {\bf 4} (1952), 100--101.

\bibitem{iwasawatamagawa1951b} K.~Iwasawa and T.~Tamagawa,   Correction: On
the paper ``On the group of automorphisms of a function field''.
 {J. Math. Soc. Japan}  {\bf 4} (1952),
203--204.


\bibitem{kantor-o'nan-seitz1972} W.M.~Kantor, M.~O'Nan and G.M.~Seitz,
 $2$-transitive groups in which the stabilizer of two points
is cyclic. {J. Algebra} {\bf 21} (1972), 17--50.

\bibitem{lehr-matignon2005} C.~Lehr and M.~Matignon, Automorphism groups
for $p$-cyclic covers of the affine line. {Compositio Math}.
\textbf{141} (2005), 1213--1237.


\bibitem{madanrosen1992} M.~Madan and M.~Rosen, The automorphism
group of a function field. {Proc. Amer. Math. Soc.} {\bf 115}
(1992), 923-929.


\bibitem{maddenevalentini1983} D.J.~Madden and R.C.~Valentini,
The group of automorphisms of algebraic function fields. {J.
Reine Angew. Math.} {\bf 343} (1983), 162--168.

\bibitem{mihailescu2004} P.~Mihauilescu, Primary cyclotomic units
and a proof of Catalan's conjecture. J. Reine Angew. Math. {\bf
572} (2004), 167--195.


\bibitem{nakajima1987}
S.~Nakajima,  $p$-ranks and automorphism groups of algebraic
curves. {Trans. Amer. Math. Soc.} {\bf 303} (1987), 595--607.


\bibitem{nakajima1987bis} S.~Nakajima, On abelian automorphism groups of
algebraic curves. {J. London Math. Soc.} {\bf 36} (1987),
23--32.

\bibitem{nakajima2002} S.~Nakajima, On automorphism groups of algebraic
curves, \emph{Current Trends in Number Theory}, Hindustan Book
Agency, New Delhi, 2002, 129--134.


\bibitem{roquette1952} P.~Roquette, \"Uber die Automorphismengruppe eines
algebraischen Funktionenk\"orpers. {Arch. Math.} {\bf 3}
(1952), 343--350.

\bibitem{roquette1970} P.~Roquette,
Absch{\"a}tzung der Automorphismenanzahl von Funktionen\-k{\"o}rpern
bei Primzahlcharakteristik. {Math. Z.} {\bf 117} (1970),
157--163.



\bibitem{schmid1938} H.I.~Schmid,
\"Uber
Automorphismen eines algebraische Funktionenk\"orpern von
Primzahlcharakteristic. {J. Reine Angew. Math.} \textbf{179}
(1938), 5-15.


\bibitem{stichtenoth1973I} H.~Stichtenoth,
\"Uber die Automorphismengruppe eines algebraischen
Funktionenk{\"o}rpers von Primzahlcharakteristik. I. Eine
Absch{\"a}tzung der Ordnung der Automorphismengruppe. {Arch.
Math.} {\bf 24} (1973), 527--544.

\bibitem{stichtenoth1973II} H.~Stichtenoth,
\"Uber die Automorphismengruppe eines algebraischen
Funktionenk{\"o}rpers von Primzahlcharakteristik. II. Ein spezieller
Typ von Funktionenk{\"o}rpern. {Arch. Math.} {\bf 24} (1973),
615--631.

\bibitem{stichtenoth1984} H.~Stichtenoth,
Zur Realisierbarkeit endlicher Gruppen als Automorphismengruppen
algebraischer Funktionenk{\"o}rper. {Math. Z.} {\bf 187}
(1984), 221--225.




\bibitem{tits1960} J.~Tits,   Les groupes simples de Suzuki et de Ree.
{S{\'e}minaire Bourbaki} {\bf 6},  Soc. Math. France, Paris,
1995, 65--82.

\bibitem{tits1962} J.~Tits,  Ovoides et groupes de Suzuki. {Arch.
Math.} {\bf 13} (1962), 187--198.

\bibitem{maddenevalentini1982} R.C.~Valentini and M.L.~Madan,  A
Hauptsatz of L.E. Dickson and Artin--Schreier extensions. {J.
Reine Angew. Math.} {\bf 318} (1980), 156--177.

\bibitem{zhang1991} J.P.~Zhang,
On the finite group with a TI Sylow $p$-subgroup. {Chinese Ann.
Math. Ser. B} {\bf 12} (1991), 147--151.
 \end{thebibliography}
\end{document}